\documentclass[12pt,twoside]{article}
\usepackage{amsmath,amssymb,amsthm}
\usepackage{enumitem}

\newtheorem{lem}{Lemma}
\newtheorem{thm}{Theorem}
\newtheorem{cor}{Corollary}
\newtheorem{prop}{Proposition}
\newtheorem{rem}{Remark}
\newtheorem{exa}{Example}

\begin{document}

\pagestyle{empty}

\title{Artin-style characterizations of the multiple gamma and sine functions}
\author{Genki Shibukawa (Kitami Institute of Technology)}
\date{\empty}
\pagestyle{plain}

\maketitle

\begin{abstract}
We extend Artin's characterizations of the gamma function to the multiple setting by means of the Kubert identity. 
We obtain uniqueness theorems for positive continuous solutions of a difference equation $G(z + 1) = \alpha(z)G(z)$ that also satisfy an Artin-type multiplication formula, remaining within the class of continuous functions under a modest initial condition. 
The results apply in particular to the Barnes multiple gamma functions and the Kurokawa multiple sine functions. 
\end{abstract}


\section{Introduction}

The gamma function
$$
\Gamma(z) := \int_{0}^{\infty} u^{z - 1} e^{-u}\, du, \qquad \mathrm{Re}(z) > 0
$$
is well known to be characterized by the Bohr--Mollerup theorem. 
This states that the gamma function is uniquely characterized by the functional equation
\begin{equation}
f(z + 1) = z f(z) \tag{FE}
\end{equation}
together with the initial condition $f(1) = 1$ and the convexity of $\log {f(z)}$ for $z > 0$ (logarithmic convexity). 
Even if one assumes continuity or continuous differentiability in addition to (FE) and the initial condition, uniqueness does not hold; the striking point is that uniqueness is restored once logarithmic convexity is imposed. 
This has been extended almost in parallel to the $q$-gamma function
$$
\Gamma_{q}(z) := (1 - q)^{1 - z} \prod_{n = 1}^{\infty} \frac{1 - q^{1 + n}}{1 - q^{z + n}}, \qquad |q| < 1
$$
by \cite{As}.

Analogous results are known for multiple gamma functions. \cite{V} showed that the multiple gamma function is characterized by
\begin{enumerate}
\item $G_{r}(z + 1) = G_{r - 1}(z) G_{r}(z)$,
\item $G_{r}(1) = 1$,
\item $\displaystyle \frac{d^{r + 1}}{dz^{r + 1}} \log G_{r}(z + 1) \geq 0$ for $z \geq 0$,
\item $G_{0}(z) = z$.
\end{enumerate}
A $q$-analogue (a multiple version of Askey's result) was given by \cite{N}. 
However, unlike the $r = 1$ case of Bohr--Mollerup and Askey, mere logarithmic convexity is insufficient in the general $r$ case; higher-order differentiability is essentially required.

On the other hand, \cite{Ar} gave several characterizations of the gamma function different from Bohr--Mollerup. 
The key is to combine the functional equation $f(z + 1) = z f(z)$ with the multiplication formula
\begin{equation}
f(z) = \frac{N^{\frac{z}{N} - \frac{1}{2}}}{(2 \pi)^{\frac{N - 1}{2}}} \prod_{j = 0}^{N - 1} f \Bigl( \frac{z + j}{N} \Bigr) \tag{A\(_N\)}
\end{equation}
satisfied by the gamma function. 
This yields, for example, the following characterizations.

\cite{Ar} Theorem 6.1. The gamma function is the unique twice continuously differentiable function that is positive for positive real $x$ and satisfies (FE) and (A\(_N\)) for $N = 2$ (written (A\(_2\))).

\cite{Ar} Theorem 6.2. The gamma function is the unique once continuously differentiable function that is positive for positive real $x$ and satisfies (FE) and (A\(_N\)) for some $N \geq 2$. 

\cite{Ar} Theorem 6.3. The gamma function is the unique continuous function that is positive for positive real $x$ and satisfies (FE) and (A\(_N\)) for every $N \geq 2$.

To prove these three theorems, Artin established the following three lemmas for the logarithmic version of (A\(_N\)):
\begin{equation}
g(z) = \sum_{j = 0}^{N - 1} g \Bigl( \frac{z + j}{N} \Bigr). \tag{lA\(_N\)}
\end{equation}

\cite{Ar} Lemma 1. Let $g$ be a twice continuously differentiable period-1 function satisfying (lA\(_N\)) for $N = 2$ (written (lA\(_2\))). Then $g \equiv 0$.

\cite{Ar} Lemma 2. Let $g$ be a once continuously differentiable period-1 function satisfying (lA\(_N\)) for some $N \geq 2$. Then $g \equiv 0$.

\cite{Ar} Lemma 3. Let $g$ be a continuous period-1 function satisfying (lA\(_N\)) for every $N \geq 2$. Then $g \equiv 0$.

As one proceeds from top to bottom, the regularity assumption on $f$ weakens from $C^{2}$ to $C^{0}$, while the multiplication formula (A\(_N\)) becomes stricter. 
In particular, Theorem 6.3 (or Lemma 3) weakens the condition on the function itself to mere continuity at the cost of requiring the multiplication formula for all integers $N \geq 2$. 
Compared with the Bohr--Mollerup theorem, which is controlled only by (FE), an initial condition, and logarithmic convexity, this may appear to be a rather complicated characterization. 

However, the situation changes in the multiple setting: the Artin-style approach (especially Theorem 6.3) is superior to the Bohr--Mollerup style in the following respects. 

1. In the multiple case, the Bohr--Mollerup approach forces logarithmic convexity into a rather strong constraint of higher-order (logarithmic) differentiability and its positivity. 
The Artin-style approach can preserve the requirement of continuous functions even in the multiple setting by adding only one extra initial condition (one may choose either higher-order differentiability or an initial condition).

2. The same characterization works not only for multiple gamma functions but also for other types of multiple special functions such as multiple sine functions.

In this paper we give a multiple analogue of Artin's characterizations of the gamma function that use the multiplication formula (A\(_N\)), and in particular the more fundamental Lemmas \ref{thm:key lem1}, \ref{thm:key lem2}, and \ref{thm:key lem3}, by means of the Kubert identity for functions of period 1:
\begin{equation}
g(z) = \sum_{j_{1}, \dots, j_{r} = 0}^{N - 1} g \Bigl( \frac{z + j_{1} + \cdots + j_{r}}{N} \Bigr) = N^{r - 1} \sum_{j = 0}^{N - 1} g \Bigl( \frac{z + j}{N} \Bigr). \tag{K\(_N\)}
\end{equation}
To this end, we give the following key lemmas. 
In what follows, unless otherwise stated, we assume that $r$ is an integer greater than or equal to 2.

\begin{lem}
\label{thm:key lem1}
\renewcommand{\thelem}{1.1}
Let $g$ be a $C^{r + 1}$ function of period 1 that satisfies the functional equation $(\mathrm{K}_2)$. Then $g \equiv 0$.
\end{lem}

\begin{lem}
\label{thm:key lem2}
\renewcommand{\thelem}{1.2}
Let $g$ be a $C^{r}$ function of period 1 that satisfies the functional equation $(\mathrm{K}_N)$ for some integer $N \geq 2$. 
Then $g \equiv 0$.
\end{lem}

\begin{lem}[\cite{M}]
\label{thm:key lem3}
\renewcommand{\thelem}{1.3}
Let $r$ be an integer greater than or equal to 2. 
Let $g$ be a continuous function of period 1 that satisfies the functional equation $(\mathrm{K}_N)$ for every integer $N \geq 2$. 
Then there exist complex constants $a, b$ such that
\begin{align}
\label{eq:g Fourier expan}
g(z)
   =
   a \sum_{n = 1}^{\infty} \frac{\sin (2 \pi n z)}{n^{r}}
   + b \sum_{n = 1}^{\infty} \frac{\cos (2 \pi n z)}{n^{r}}.
\end{align}
\end{lem}

Lemma~\ref{thm:key lem3} is due to \cite{M}. 
Its proof is a well-known consequence of Fourier series expansion, but we revisit it in this paper in view of the following points.
\begin{itemize}
\item There are several variations of the additional conditions that force $g \equiv 0$ from this representation.
\item Surprisingly, no mention seems to have been made of its application to the characterization of multiple gamma functions or multiple sine functions.
\end{itemize}

As an application we state characterizations of a class of multiple special functions that includes the multiple gamma functions and the multiple sine functions. 
Let $G(z)$ be a continuous function that is positive on the positive reals and satisfies the functional equation
\begin{equation}
G(z + 1) = \alpha(z) G(z). \tag{GE}
\end{equation}
Assume further that, for every integer $N \geq 2$, $G$ satisfies a multiplication formula of the form
\begin{equation}
G(z) = \beta(z, N) \prod_{j_{1}, \dots, j_{r} = 0}^{N - 1} G \Bigl( \frac{z + j_{1} + \cdots + j_{r}}{N} \Bigr). \tag{\(B_N\)}
\end{equation}

\begin{thm}
\label{thm:multiple analogue of A1}
If functions $G_{1}$, $G_{2} \in C^{r + 1}(\mathbb{R}_{>0})$ satisfy (GE) and the multiplication formula ($B_2$) for $N = 2$, then $G_{1} = G_{2}$. 
\end{thm}

\begin{thm}
\label{thm:multiple analogue of A2}
If functions $G_{1}$, $G_{2} \in C^{r}(\mathbb{R}_{>0})$ satisfy (GE) and the multiplication formula ($B_N$) for some integer $N \geq 2$, then $G_{1} = G_{2}$. 
\end{thm}

\begin{thm}
\label{thm:multiple analogue of A3}
\begin{enumerate}
\item[(1)] If continuous functions $G_{1}$, $G_{2} \in C(\mathbb{R}_{>0})$ satisfy (GE) and the multiplication formula ($B_N$) for every integer $N \geq 2$, then there exist complex constants $a$  and $b$ such that
$$
G_{1}(z) = \exp \Biggl( a \sum_{n = 1}^{\infty} \frac{\sin (2 \pi n z)}{n^r} + b \sum_{n = 1}^{\infty} \frac{\cos (2 \pi n z)}{n^r} \Biggr)G_{2}(z).
$$
\item[(2)] If continuous functions $G_{1}$, $G_{2} \in C(\mathbb{R}_{>0})$ satisfy (GE), the multiplication formula ($B_N$) for every integer $N \geq 2$, and $G_{1}(1) = G_{2}(1)$ and $G_{1}(1/m) = G_{2}(1/m)$ for some integer $m \geq 3$, then $G_{1}(z) = G_{2}(z)$.
\end{enumerate}
\end{thm}

The paper is organized as follows. 
Section~2 introduces the Barnes-type multiple zeta functions, multiple gamma functions, and multiple sine functions and recalls their basic properties. 
Section~3 proves or recalls Lemmas~\ref{thm:key lem1},~\ref{thm:key lem2}, and~\ref{thm:key lem3}, which are the multiple versions of Artin's Lemmas 1, 2, and 3, and states some corollaries. 
Finally, we prove Theorems~\ref{thm:multiple analogue of A1}, ~\ref{thm:multiple analogue of A2} and~\ref{thm:multiple analogue of A3}, and discusses characterizations of more general multiple special functions that include the multiple gamma functions and multiple sine functions.

\section{Preliminaries}
Refer to \cite{B} and \cite{KK} for further details in this section. 
For any complex number $a$ we take the argument in the range $-\pi < \arg a \leq \pi$. 
For a positive integer $r$, let $\omega = (\omega_{1}, \dots, \omega_{r}) \in (\mathbb{C}^{\times})^{r}$ be a tuple of non-zero complex parameters. 
The Barnes-type multiple zeta function is defined by
$$
\zeta_{r}(s, z \mid \omega) := \sum_{m_{1}, \dots, m_{r} \geq 0} \frac{1}{(m_{1} \omega_{1} + \cdots + m_{r} \omega_{r} + z)^s}, \qquad \mathrm{Re}(s) > r,
$$
where $z, \omega_{1}, \dots, \omega_{r}$ are assumed to satisfy
$$
\max \{ \arg z, \arg \omega_{1}, \dots, \arg \omega_{r} \} - \min \{ \arg z, \arg \omega_{1}, \dots, \arg \omega_{r} \} < \pi
$$
(throughout the paper, unless otherwise stated, we always assume this condition).

It is well known that $\zeta_{r}(s, z \mid \omega)$ admits a meromorphic continuation to the whole $s$-plane, with poles only at $s = 1, 2, \dots, r$, each of which is simple. 
In particular, it is holomorphic at $s = 0$, and we define the Barnes multiple gamma function by
$$
\Gamma_{r}(z \mid \omega) := \exp \Biggl( \Biggl. \frac{\partial}{\partial s} \zeta_{r}(s, z \mid \omega) \Biggr|_{s = 0} \Biggr).
$$
When $r = 1$ one has
$$
\Gamma_{1}(z \mid \omega_1) = \frac{\omega_{1}^{\frac{z}{\omega_{1}} - \frac{1}{2}}}{\sqrt{2 \pi}}\Gamma\left(\frac{z}{\omega_{1}}\right).
$$

The difference relation for the multiple zeta function
$$
\zeta_{r}(s, z + \omega_{i} \mid \omega) = \zeta_{r}(s, z \mid \omega) - \zeta_{r - 1}(s, z \mid \widehat{\omega}^{(i)}), \qquad i = 1, \dots, r,
$$
where $\widehat{\omega}^{(i)} = (\omega_{1}, \dots, \omega_{i - 1}, \omega_{i + 1}, \dots, \omega_{r}) \in (\mathbb{C}^\times)^{r - 1}$, implies the difference relation
$$
\Gamma_{r}(z + \omega_{i} \mid \omega) = \Gamma_{r - 1}(z \mid \widehat{\omega}^{(i)})^{-1} \Gamma_{r}(z \mid \omega).
$$

Using the multiple gamma function we define Kurokawa's multiple sine function by
$$
S_{r}(z \mid \omega) := \Gamma_{r}(z \mid \omega)^{-1} \Gamma_{r}(\omega_{1} + \cdots + \omega_{r} - z \mid \omega)^{(-1)^{r}}.
$$
It satisfies $S_{1}(z \mid \omega_1) = 2 \sin (\pi z / \omega_{1})$ and the difference relation
$$
S_{r}(z + \omega_{i} \mid \omega) = S_{r - 1}(z \mid \widehat{\omega}^{(i)})^{-1} S_{r}(z \mid \omega).
$$

The following multiplication formulas, which are the multiple analogues of Artin's formula (A\(_N\)), are known:
\begin{align*}
\Gamma_{r}(z \mid \omega) & = N^{- \zeta_{r}(0, z \mid \omega)} \prod_{j_{1}, \dots, j_{r} = 0}^{N - 1} \Gamma_{r} \Biggl( \frac{z + j_{1} \omega_{1} + \cdots + j_{r} \omega_{r}}{N} \Biggm| \omega \Biggr), \\
S_{r}(z \mid \omega) & = \prod_{j_{1}, \dots, j_{r} = 0}^{N - 1} S_{r} \Biggl( \frac{z + j_{1} \omega_{1} + \cdots + j_{r} \omega_{r}}{N} \Biggm| \omega \Biggr).
\end{align*}

When $\omega_{1} = \cdots = \omega_{r} = 1$ we omit the parameter $\omega$ and simply write
$$
\zeta_{r}(s, z), \qquad \Gamma_{r}(z), \qquad S_{r}(z).
$$

\begin{rem}
\label{rem:normalizations}
(1) An alternative normalization
$$
G_{r}(z \mid \omega) := \exp \Biggl( (-1)^{r - 1} \Biggl. \frac{\partial}{\partial s} \zeta_{r}(s, z \mid \omega) \Biggr|_{s = 0} \Biggr)
$$
is sometimes used. In this normalization the difference equation becomes
$$
G_{r}(z + \omega_{i} \mid \omega) = G_{r}(z \mid \omega) G_{r - 1}(z \mid \widehat{\omega}^{(i)}).
$$

(2) The value $\zeta_{r}(0, z \mid \omega)$ appearing in the multiplication formula is given in terms of the multiple Bernoulli polynomials defined by the generating function
$$
\frac{u^{r} e^{z u}}{(e^{\omega_{1} u} - 1) \cdots (e^{\omega_{r} u} - 1)} = \sum_{n = 0}^{\infty} B_{r, n}(z \mid \omega) \frac{u^{n}}{n!}
$$
via the relation
$$
\zeta_{r}(0, z \mid \omega) = (-1)^{r} \frac{B_{r, r}(z \mid \omega)}{r!}.
$$
In particular it is a polynomial of degree $r$ in $z$.

(3) When $r \geq 3$ and $\omega_{i} / \omega_{j} \notin \mathbb{Q}$ for all $i \neq j$, Jacobi's theorem implies that there are no non-constant meromorphic functions with three or more independent periods. 
Then, any solution of the system of difference equations
$$
f_{r}(z + \omega_{j} \mid \omega) = f_{r - 1}(z \mid \widehat{\omega}^{(j)}) f_{r}(z \mid \omega), \qquad j = 1, \dots, r
$$
is unique up to a constant factor. 
Consequently, for $r \geq 3$ the uniqueness problem for multiple special functions becomes non-trivial only when the parameters are rationally dependent. 
In the present paper we treat the most degenerate case $\omega_{1} = \cdots = \omega_{r} = 1$.
\end{rem}

Finally, although unrelated to the multiple zeta function, we prove the following proposition, which will be needed below.

\begin{prop}
\label{thm:non zero prop}
For any real number $r > 1$ and any integer $N \geq 3$,
$$
\sum_{n = 1}^{\infty} \frac{\sin \Bigl( \frac{2 \pi n}{N} \Bigr)}{n^{r}} > 0.
$$
In particular the sum is non-zero.
\end{prop}

\begin{proof}
Under the stated hypotheses the series converges absolutely, and we may rearrange it as
\begin{align*}
\sum_{n = 1}^{\infty} \frac{\sin \Bigl( \frac{2 \pi n}{N} \Bigr)}{n^{r}}
& = \sum_{m = 0}^{\infty} \sum_{j = 1}^{N} \frac{\sin \Bigl( \frac{2 \pi j}{N} \Bigr)}{(m N + j)^{r}} \\
& = \sum_{j = 1}^{\lfloor \frac{N}{2} \rfloor} \sin \Bigl( \frac{2 \pi j}{N} \Bigr) \sum_{m = 0}^{\infty} \Biggl( \frac{1}{(m N + j)^{r}} - \frac{1}{(m N + N - j)^{r}} \Biggr).
\end{align*}
Each inner sum is strictly positive, which yields the claim.
\end{proof}

\section{Proofs of main results}

We prove Lemmas~\ref{thm:key lem1} and~\ref{thm:key lem2}. The arguments rely on Artin's Lemmas 1 and 2 (see~\cite{Ar}).

\begin{proof}[Proof of Lemma~\ref{thm:key lem1}]
Assume that $g$ is a $C^{r + 1}$ function of period 1 satisfying $(\mathrm{K}_2)$. 
Differentiating both sides $r - 1$ times with respect to $z$ yields
\begin{equation}
g^{(r - 1)}(z) = g^{(r - 1)} \Bigl( \frac{z}{2} \Bigr) + g^{(r - 1)} \Bigl( \frac{z + 1}{2} \Bigr). \tag{\(\mathrm{K}_2^{(r - 1)}\)}
\end{equation}
Here $g^{(r - 1)}(z)$ is a twice continuously differentiable period 1 function on $[0, 1]$, so Artin's Lemma 1 implies $g^{(r - 1)}(z) \equiv 0$.

Integrating $r - 1$ times shows that $g(z)$ is a polynomial of degree at most $r - 1$. 
Since it is periodic of period 1, it must be constant: $g(z) = c$. 
Substituting into $(\mathrm{K}_2)$ gives $c = 2^{r} c$. 
Since $r \geq 2$, we conclude that $c = 0$.
\end{proof}

\begin{proof}[Proof of Lemma~\ref{thm:key lem2}]
Assume that $g$ is a $C^{r}$ function of period 1 satisfying $(\mathrm{K}_{N})$ for some integer $N \geq 2$. 
Differentiating both sides $r - 1$ times with respect to $z$ yields
\begin{equation}
g^{(r - 1)}(z) = \sum_{j = 0}^{N - 1} g^{(r - 1)} \Bigl( \frac{z + j}{N} \Bigr). \tag{\(\mathrm{K}_N^{(r - 1)}\)}
\end{equation}
Here $g^{(r - 1)}(z)$ is a $C^{1}$ function of period 1, so Artin's Lemma 2 implies $g^{(r - 1)} \equiv 0$. 
The same argument as in the previous proof then shows that $g \equiv 0$.
\end{proof}

Lemma~\ref{thm:key lem3} is Lemma 1 of~\cite{M}, so we do not reproduce the proof. 
We record, however, several natural initial conditions that force $g \equiv 0$.

\begin{cor}
\label{thm:cor of key lem}
Let $g$ be a continuous function of period 1 satisfying $(\mathrm{K}_N)$ for every integer $N \geq 2$.
\begin{enumerate}
\item[(1)] If $g$ is of class $C^{1}$ and $g(0) = g'(0) = 0$, then $g \equiv 0$.
\item[(2)] If $g(0) = 0$ and $g(1 / m) = 0$ for some integer $m \geq 3$, then $g \equiv 0$.
\end{enumerate}
\end{cor}

\begin{proof}
Since $g$ satisfies the hypotheses of Lemma~\ref{thm:key lem3}, there exist complex constants $a, b$ such that~\eqref{eq:g Fourier expan} holds.

(1) The condition $g(0) = 0$ together with $\zeta(r) > 1$ forces $b = 0$, so that
\begin{equation}
\label{eq:g odd Fourier}
g(z) = a \sum_{n = 1}^{\infty} \frac{\sin (2 \pi n z)}{n^r}.
\end{equation}
Differentiating term by term gives
$$
g'(z) = 2 \pi a \sum_{n = 1}^{\infty} \frac{\cos (2 \pi n z)}{n^{r - 1}}.
$$
If $r = 2$, the right-hand side is discontinuous at the integers whenever $a \neq 0$, contradicting the $C^1$ assumption, so $a = 0$. 
If $r > 2$, the condition $g'(0) = 0$ together with $\zeta(r - 1) > 1$ again forces $a = 0$. 

(2) From $g(0) = 0$ we obtain~\eqref{eq:g odd Fourier} as before. 
Substituting $z = 1 / m$ and applying Proposition~\ref{thm:non zero prop} yields $a = 0$.
\end{proof}

\begin{rem}
In the multiple case, Milnor's result (Lemma~\ref{thm:key lem3}) requires additional initial conditions to conclude that $g \equiv 0$. 
By contrast, Artin's Lemma 3 (the classical $r = 1$ case) needs no such conditions. 
The reason is that the Fourier series
$$
a \sum_{n = 1}^{\infty} \frac{\sin (2 \pi n z)}{n} + b \sum_{n = 1}^{\infty} \frac{\cos (2 \pi n z)}{n}
$$
fails to converge absolutely; the resulting discontinuities force the coefficients to vanish under the mere assumption of continuity.
\end{rem}

From Lemmas~\ref{thm:key lem1},~\ref{thm:key lem2}, and~\ref{thm:key lem3}, we prove Theorems~\ref{thm:multiple analogue of A1}, ~\ref{thm:multiple analogue of A2} and~\ref{thm:multiple analogue of A3}.

\begin{proof}[Proof of Theorems~\ref{thm:multiple analogue of A1}, ~\ref{thm:multiple analogue of A2} and~\ref{thm:multiple analogue of A3}]
For $z > 0$ set
$$
g(z) := \log \Biggl( \frac{G_1(z)}{G_2(z)} \Biggr).
$$
Then $g$ is a continuous function of period 1. 
The multiplication formula ($B_N$) implies the Kubert identity ($\mathrm{K}_N$). 
The claims now follow from Lemmas~\ref{thm:key lem1},~\ref{thm:key lem2},~\ref{thm:key lem3} and Corollary~\ref{thm:cor of key lem}.
\end{proof}

\begin{exa}
Taking $\alpha(z) = \Gamma_{r - 1}(z)^{-1}$ and $\beta(z, N) = N^{- \zeta_r(0, z)}$, the Barnes multiple gamma function $\Gamma_r(z)$ is a continuous positive solution on $\mathbb{R}_{> 0}$ of both (GE) and ($B_N$) for every $N \geq 2$. 
Consequently, by Theorem~\ref{thm:multiple analogue of A3}, if $G(z)$ satisfies (GE) and ($B_N$) for every integer $N \geq 2$, is continuous, and satisfies $G(1) = \Gamma_r(1)$ and $G(1 / m) = \Gamma_r(1 / m)$ for some integer $m \geq 3$, then $G(z) = \Gamma_r(z)$.

The same argument with $\alpha(z) = S_{r - 1}(z)^{-1}$ and $\beta(z, N) = 1$ shows that $G(z) = S_r(z)$.
\end{exa}

\begin{rem}
More generally, the same characterizations hold for the Milnor multiple gamma functions
$$
\Gamma_{r, k}(z) := \exp \Biggl( \Biggl. \frac{\partial}{\partial s} \zeta_r(s, z) \Biggr|_{s = -k} \Biggr),
$$
provided $r - k > 1$. 
In this case the undetermined factor appearing in Theorem~\ref{thm:multiple analogue of A3} is replaced by
$$
\exp \Biggl( a \sum_{n = 1}^{\infty} \frac{\sin (2 \pi n z)}{n^{r - k}} + b \sum_{n = 1}^{\infty} \frac{\cos (2 \pi n z)}{n^{r - k}} \Biggr).
$$
\end{rem}

It remains an open question whether the Artin-style characterizations via the Kubert identity developed in this paper admit $q$-analogue, like the $q$-Bohr--Mollerup theorems of \cite{As} and \cite{N}. 
Another possible direction is to weaken the regularity assumptions appearing in Lemmas~\ref{thm:key lem1} and~\ref{thm:key lem2}, in view of the optimality results obtained by \cite{K} for the classical case. 
We leave both directions for future work.

\section*{Acknowledgment}
This work is supported by JSPS KAKENHI Grant Number 26K06837. 
English language editing was assisted by Grok, an AI developed by xAI. 


\bibliographystyle{amsplain}

\medskip
\begin{flushleft}
Genki Shibukawa\\ 
Kitami Institute of Technology\\
165, Koen-cho, Kitami, Hokkaido 090-8507\\
Japan\\
E-mail: g-shibukawa@mail.kitami-it.ac.jp
\end{flushleft}

\end{document}